\documentclass[11pt,reqno]{amsart}
\usepackage{fullpage}

\usepackage{graphicx}
\usepackage[draft]{hyperref}
\usepackage{amsmath,amsopn,amssymb,amsfonts,stmaryrd}
\usepackage{verbatim}
\usepackage{amsthm}
\usepackage{mathtools}
\usepackage{color}
\usepackage{enumitem}
\usepackage[framemethod=TikZ]{mdframed}
\usepackage{bbm}
\usepackage{mathrsfs}
\usepackage{booktabs}
\usepackage{caption}
\usepackage{bm}

\theoremstyle{plain}
\newtheorem{theorem}{Theorem}[section]

\newtheorem{lemma}[theorem]{Lemma}
\newtheorem{proposition}[theorem]{Proposition}
\theoremstyle{definition}

\newtheorem*{remark}{Remark}

\theoremstyle{remark}

\newcommand{\IR}{{\mathbb R}}
\newcommand{\IC}{{\mathbb C}}
\newcommand{\N}{{\mathbb N}}

\newcommand{\SU}{\mathrm{SU}}

\numberwithin{equation}{section}

\allowdisplaybreaks

\newcommand{\mat}[1]{\left( \begin{matrix} #1 \end{matrix} \right)}

\def\lp{\left(}
\def\rp{\right)}

\def\a{\alpha}
\def\b{\beta}
\def\d{\delta}

\def\z{\zeta}

\def\s{\sigma}

\setlist[itemize]{noitemsep, topsep=0pt}

\makeatletter
\newcommand{\vast}{\bBigg@{2}}
\newcommand{\Vast}{\bBigg@{5}}
\makeatother

\title{An asymptotic formula for the number of $n$-dimensional representations of $\SU(3)$}
\author{Kathrin Bringmann and Johann Franke}
\address{University of Cologne, Department of Mathematics and Computer Science, Weyertal 86-90, 50931 Cologne, Germany}
\email{kbringma@math.uni-koeln.de }
\email{{jfrank12@uni-koeln.de}}

\begin{document}
\maketitle

\section{Introduction and statement of results}

The special unitary group $\SU(2)$ has (up to equivalence) one irreducible representation $V_k$ of each dimension $k\in \N$. Each $n$-dimensional representation $\bigoplus_{k=1}^\infty r_k V_k$ corresponds to a unique partition 
\begin{align}
\label{partition} n = \lambda_1 + \lambda_2 + \cdots + \lambda_r, \qquad 1 \leq \lambda_1 \leq \lambda_2 \leq \ldots \leq \lambda_r \leq n,
\end{align}
such that $r_k$ denotes the number of $k$'s in \eqref{partition}. As a result, the number of representations equals  $p(n)$, the number of integer partitions of $n$. The partition function has no elementary closed formula, nor does it satisfy any finite order recurrence. However, with $p(0) := 1$, its generating function has the following classical product expansion
\begin{align}
\label{partitionproduct} \sum_{n=0}^\infty p(n)q^n = \prod_{n=1}^\infty \frac{1}{1-q^n}.
\end{align}
In \cite{HardyRama}, Hardy and Ramanujan showed the asymptotic formula 
\begin{align}
\label{partitionasy} p(n) \sim \frac{1}{4\sqrt{3} n} \exp\left( \pi \sqrt{\frac{2n}{3}}\right), \qquad n \rightarrow \infty.
\end{align}
For their proof they introduced the now so-called Hardy--Ramanujan Circle Method, that uses modular type transformations to obtain a divergent asymptotic expansion whose truncations approximate $p(n)$ up to small errors. A later refinement by Rademacher \cite{Rademacher} provides a convergent series for $p(n)$. 

It is natural to ask whether the above correspondence between the number of representations of a group and a  partition function can be generalized. The obvious next case is the unitary group $\SU(3)$, whose irreducible representations are a family of representations $W_{j,k}$ indexed by pairs of positive integers. Note that (see Chapter 5 of \cite{Hall}) that $\dim (W_{j,k}) = \tfrac12 jk(j+k)$. Like in the case of $\SU(2)$, a general $n$-dimensional representation decomposes into a sum of these $W_{j,k}$, again each with some multiplicity. So analogous to \eqref{partitionproduct}, it is easy to see that the numbers $r(n)$ of $n$-dimensional representations, again with $r(0) := 1$, have the generating function 
\begin{equation}\label{G}
 G(q):= \sum_{n=0}^{\infty} r(n)q^n=\prod_{j,k\ge1} \frac{1}{1-q^{\frac{jk(j+k)}{2}}} = 1 + q + q^2 + 3q^3 + 3q^4 + 3q^5 + 8q^6 + 8q^7 + \cdots .
\end{equation}
In \cite{Ro}, Romik proved the following analogon of formula \eqref{partitionasy} for the sequence $r(n)$: 
\begin{equation} \label{rn}
	r(n) \sim \frac{C_0}{n^{\frac 35}} \exp\left(A_1n^{\frac 25}-A_2n^{\frac{3}{10}}-A_3n^{\frac 15}- A_4n^{\frac{1}{10}}\right) \quad (n \to \infty),
\end{equation}
for $A_1,A_2,A_3,A_4,C_0\in\IR$, that are defined in the notation section. 
Romik asked for lower order terms in the asymptotic expansion of $r(n)$. We answer his question in the following theorem.

\begin{theorem}\label{Main}
	Let $L \in\N_0$. We have, as $n \rightarrow \infty$,
	\begin{align*}
	r(n) = \frac{1}{n^{\frac 35}} \left(\sum_{j=0}^{L}\frac{C_j}{n^\frac{j}{10}} + O_L\left( n^{-\frac{L}{10}-\frac{3}{80}}\right)\right) \exp\left(A_1n^{\frac 25}-A_2n^{\frac{3}{10}}-A_3n^{\frac 15}- A_4n^{\frac{1}{10}}\right),
	\end{align*}
	where the constants $C_j$ do not depend on $L$ and $n$ and can all be calculated explicitly. 
\end{theorem} 

\begin{remark} In \eqref{C1} and \eqref{C2} we compute some of these constants explicitly.
\end{remark} 

In order to prove \eqref{rn} and Theorem \ref{Main}, one has to carefully study the related Dirichlet series 
\begin{equation}\label{omega}
\omega(s) := \sum_{j,k \ge 1} \frac{1}{(2 \dim (W_{j,k}))^{s}} = \sum_{k,j\ge 1} \frac{1}{k^s j^s (k+j)^s},
\end{equation}
that converges absolutely for all $s \in \mathbb{C}$ with $\mathrm{Re}(s) > \frac23$. It is natural to ask whether $\omega$ can be continued to a meromorphic function on the entire complex plane and satisfies a functional equation. While the second question seems out of reach, the first one was answered positivly by Matsumoto \cite{Mat}. Romik \cite{Ro} then studied the properties of $\omega$ beyond its abscissa of convergence more closely. Although $\omega(s)$ does not  seem to possess a simple functional equation, it turns out that it has ``trivial zeros" at $s \in-\N$. The key identity behind this is
\begin{align}
\label{zetaid} \zeta(6n+2) = \frac{2(4n+1)!}{(6n+1)(2n)!^2}  \sum_{k=1}^n \frac{\binom{2n}{2k-1}}{\binom{6n}{2n+2k-1}} \zeta(2n+2k) \zeta(4n-2k+2) \qquad (n \in \N),
\end{align}
where as usual the {\it Riemann $\z$-function} is
$
 \zeta(s) := \sum_{n = 1}^\infty \frac{1}{n^s} \ \left(\mathrm{Re}(s) > 1\right).
$
Romik \cite{Ro} proved that \eqref{zetaid} is equivalent to the fact that $\omega(s)$ has zeros at $s\in -\N$. Both the trivial zeros and the distribution of poles of $\omega$ are of importance in the proof of Theorem \ref{Main}. A key tool for the proof of the main theorem is to give uniform estimates for the function $\omega$ on vertical lines.

The paper is organized as follows. After recalling some  preliminaries in Section \ref{sec:pre}, we prove in Section 3 an asymptotic expansion of the function $G(q)$ near $q=1$ with an explicit formula for the occuring error term. In Section 4 we give a unifom upper bound for this error which is below used to show that it is negligible. Also, we prepare the error analysis for the error tails of the major arc integral when using the Saddle Point Method. In Section 5 we give a power series expansion for the Saddle Points and show that they depend on the parameter $n$, but stay bounded as $n$ increases (in fact, they even converge to 1). This saddle point function is then used in Section 6 to investigate the asymptotic expansion of $G(q)$ near to $q=1$ in further detail. Finally, in Section 7 we use Wright's Circle Method and the Saddle Point Method to prove Theorem \ref{Main}. While the major arc integral can be evaluated using the preliminary sections, we use a lemma by Romik \cite{Ro} to deal with the minor arcs. In Section 8 we finally discuss some open problems and related questions. 

\section*{Acknowledgements}
The authors thank Caner Nazaroglu for  checking our numerical calculations. The first author has received funding from the European Research Council (ERC) under the European Union’s Horizon 2020 research and innovation programme (grant agreement No. 101001179), and the authors are partially supported by the Alfried Krupp prize. 

\section*{Notation}

In this section, we provide frequently used notation. For $\b\in\IR$, we denote by  $\{\b\}:=\b-\lfloor\b\rfloor$ the \textit{fractional part} of $\b$. As usual, $\mathbb{H}:=\{ \tau \in \IC : \mathrm{Im}(\tau) > 0\}$. For  $0 < \delta \leq \frac{\pi}{2}$, we define the cone
	\begin{align*}
	\mathcal{C}_\delta := \left\{z \in -i \mathbb{H} : |\mathrm{Arg}(z)| \leq \frac{\pi}{2} - \delta \right\},
	\end{align*}
	where $\mathrm{Arg}$ is the principal branch of the complex argument. Moreover, for $\d\in\IR$, let
	\begin{align*}
	\mathbb{C}_\delta := \{ w \in \IC : \mathrm{Im}(w) \leq 1 - \delta \}.
	\end{align*}
	We denote for $r>0$, $B_r(z):=\{w\in\IC : |w - z| < r\}$. For $z_1, z_2 \in \IC$ we define the line connecting them as 
	\begin{align*}
	[z_1, z_2] := \{ z = z_1 + \lambda (z_2-z_1) : 0 \leq \lambda \leq 1\}.
	\end{align*}
We frequently make use of the notation $f(x_1, ..., x_k) \ll g(x_1, ..., x_k)$ for $f, g: D \to \IC$, where $D \subset \IC$, which is equivalent to $\sup_{(x_1, ..., x_k) \in D} |\frac{f(x_1, ..., x_k)}{g(x_1, ..., x_k)}| < \infty$. In the case, that the constant depends on additional parameters $a,b,c,...$, we write $f(x_1, ..., x_k) \ll_{a,b,c,...} g(x_1, ..., x_k)$. We equivalently write $f(x_1 ,..., x_k) = O(g(x_1, ..., x_k))$ and $f(x_1, ..., x_k) = O_{a,b,c,...}(g(x_1, ..., x_k))$, respectively.

We also use several constants, given as
\begin{align*}
X & := \left( \tfrac19 \Gamma\left( \tfrac{1}{3}\right)^2 \zeta\left( \tfrac{5}{3}\right)\right)^{\frac{3}{10}} = 1.17117 \ldots, \  
Y  := -\sqrt{\pi} \zeta\left( \tfrac12 \right) \zeta\left( \tfrac32 \right) = 6.76190 \ldots , \\
A_1  &:= 5X^2 = 6.85826 \ldots , \
A_2  := \frac{Y}{X} = 5.7736 \ldots, \
A_3 := \frac{3Y^2}{80X^4} = 0.91134 \ldots, \\
A_4  &:= \frac{11Y^3}{3200X^7} = 0.35163 \ldots, \
C_0  := \frac{2\sqrt{3\pi}}{\sqrt{5}} X^{\frac13} \exp\left( - \frac{Y^4}{2560X^{10}} \right) = 2.44629 \ldots. 
\end{align*}
Here, $\Gamma(s)$ denotes the {\it Gamma function}, that is defined by 
\begin{align*}
\Gamma(s) := \int_0^\infty e^{-t} t^{s-1} dt \qquad (\mathrm{Re}(s) > 0).
\end{align*}

\section{Preliminaries}\label{sec:pre}

In this section we recall and prove results required for this paper. The following properties of the gamma function are well-known; the proof uses \cite{Andrews,Tenenbaum}. 

\begin{theorem} \label{GammaCollect} 
	
	\noindent
	\begin{enumerate}[leftmargin=*,label=\textnormal{(\arabic*)}]
		\item The Gamma function has a meromorphic continuation to $\IC$, with simple poles only at $s \in -\N_0$, and satisfies the function equation 
		\begin{equation*}\label{GammaFunc}
			\Gamma(s+1) = s\Gamma(s).
		\end{equation*}
		\item We have the Stirling approximation 
		\begin{align*} 
		 \Gamma(s+1) \sim \sqrt{2\pi s} \cdot s^s e^{-s}, \qquad |s| \rightarrow \infty. 
		\end{align*}
		\item Writing $s = \sigma + it$ and $\sigma \in I$ for a compact interval $I \subset [\frac12, \infty)$, we uniformly have, for $t \in \IR$,
		\begin{equation*}
		\label{Gamma}  \max\left\{1, |t|^{\sigma - \frac12}\right\} e^{- \frac{\pi |t|}{2} } \ll_I |\Gamma(\sigma + it)| \ll_I  \max\left\{1, |t|^{\sigma - \frac12}\right\}e^{- \frac{\pi|t|}{2} }.
		\end{equation*}
		The relation holds for all compact intervals $I \subset \IR$ if $|t| \ge 1$.
		\item For all $s \in \IC$ we have the identity between meromorphic functions
		\begin{equation*}\label{EulerRef}
		\Gamma(s)\Gamma(1-s) = \frac{\pi}{\sin(\pi s)}.
		\end{equation*}
		\item For $m \in \N_0$, we have 
		\begin{equation*}
		\Gamma\left(m+\frac12\right) = \frac{(2m)!}{ 4^m m!} \sqrt{\pi},\qquad
		\Gamma\left(-m-\frac12\right) = \frac{(-4)^{m+1} \left( m+1 \right)!}{(2(m+1))!} \sqrt{\pi}.
		\end{equation*}
		\item For all $x,y > -\frac12$ we have 
		\[
			 \Gamma(x+1)\Gamma(y+1) \ll \Gamma(x+y+1).
		\] 
		\end{enumerate}
	\end{theorem}

Define for $x > 0$ and $s \in \IC$ the {\it incomplete Gamma function} 
\begin{align*}
\Gamma(s;x) := \int_x^\infty t^{s-1} e^{-t} dt.
\end{align*}
It satisfies following important properties. 

\begin{theorem} \label{IncompleteGamma}
	\ \begin{enumerate}[leftmargin=*,label=\textnormal{(\arabic*)}]
					\item The function $s\mapsto\Gamma(s;x)$ defines an entire function. One has the asymptotic behavior
			\begin{equation*} 
			\label{IncomGammaAsy} \Gamma(s;x) \sim x^{s-1} e^{-x}, \qquad x \rightarrow \infty.
			\end{equation*}
					\item For $n\in\N$ and $x\in\IR$ we have the inequality 
			\begin{align*}
			\Gamma(n;x) \leq n! \max\left\{1, x^{n-1}\right\} e^{-x}.
			\end{align*}
			\end{enumerate}
	\end{theorem}

We also consider, for $s,z \in \IC$ with $s \notin -\N$, the {\it generalized Binomial coefficient} defined by 
\begin{align*}
\binom{s}{z} := \frac{\Gamma(s+1)}{\Gamma(z+1)\Gamma(s-z+1)}.
\end{align*} 
The next lemma follows by a straightforward calculation. 

\begin{lemma} \label{prep} We have the following:
	\begin{enumerate}[leftmargin=*,label=\textnormal{(\arabic*)}]
		\item We uniformly have for all $\b >0$  and $T\ge0$
		\begin{align*}
		\left| \prod\limits_{j=1}^{\lfloor \b \rfloor}(\{\b\}+j+iT)\right| \ll \Gamma(\b + 1) \max\left\{ 1, T^{\lfloor \b \rfloor }\right\}.
		\end{align*}
		\item We uniformly have, for $k \in \N_0$, $\alpha \in \N + \frac14$, and $T \geq 0$, 
		\begin{align*}
		\left|\binom{k-\a-1+iT}{k}\right| = \left| \frac{\Gamma(k-\a+iT)}{\Gamma(-\a+iT)k!} \right| \ll \left| \binom{\alpha}{k}\right| \max\left\{1, T^k\right\}.
		\end{align*}
		\item Let $0 < x < 1$.  Then there exist $m,M>0$ only dependent on $x$, such that for all $k\in\N$ 
		\begin{align*}
		mk^{x-1} \leq \binom{-x}{2k} \leq Mk^{x-1}.
		\end{align*}
	\end{enumerate}
\end{lemma}

We also need some analytic properties of the Riemann $\z$-function; see \cite{Apostol,Brue,Tenenbaum}.

\begin{theorem}\label{ZetaFunc} 
	\ \begin{enumerate}[leftmargin=*,label=\textnormal{(\arabic*)}]
	\item The $\z$-function has a meromorphic continuation to $\IC$ with only a simple pole at $s=1$ with residue $1$. For $s \in \IC$ we have (as identity between meromorphic functions)
	\begin{equation*}\label{zetafunc}
		\zeta(s) = 2^s \pi^{s-1} \sin\left( \frac{\pi s}{2}\right) \Gamma(1-s) \zeta(1-s).
	\end{equation*}   
	
	\item  Fix $c > 0, \ \s_0\in\IR$. Then we have uniformly for $|t|\ge1$ and $\s\ge\s_0$ 
	\begin{equation*}
	\label{zetadown} \zeta(\sigma + it) \ll_c |t|^{\mu(\sigma) + c},
	\end{equation*}
	where
	\begin{align*}
	\mu(\sigma) := \begin{cases} 0 & \qquad \mathrm{if } \ \sigma \geq 1, \\ \frac12 (1 - \sigma) & \qquad \mathrm{if } \ 0 \leq \sigma < 1, \\ \frac12 - \sigma & \qquad \mathrm{if } \ \sigma < 0.\end{cases}
	\end{align*}
	In particular, we obtain, uniformly in $T \ge 0$,
	\begin{align*}
	\zeta\left( \frac12 + iT \right) \ll \max\left\{ 1, T^{\frac12}\right\}.
	\end{align*}
	\end{enumerate}
\end{theorem}

To use the Saddle Point Method, we also require the following straightforward approximation. 

\begin{lemma} \label{Gauss} Let $\lambda_n$ be an increasing unbounded sequence of positive real numbers, $B > 0$, and $P$ a polynomial of degree $m\in \N_0$. Then we have 
	\begin{align*}
	\int_{-\lambda_n}^{\lambda_n} P(x) e^{-Bx^2} dx = \int_{-\infty}^{\infty} P(x) e^{-Bx^2} dx + O_{B,P}\left( \lambda_n^{\frac{m-1}{2}} e^{-B\lambda_n^2} \right).
	\end{align*}
\end{lemma}

\section{An Asymptotic expansion  of $G$}

They key to understand the asymptotic behavior of the sequence $r(n)$ is to study the function $z\mapsto G(e^{-z})$ (see \eqref{G}) next $z=0$. For this we require the following properties of $\omega$.

\begin{theorem}[Romik \cite{Ro}, Theorems 1.2 and 1.3] \label{omegapro} We have the following:
	\begin{enumerate}[leftmargin=*]
		\item[\textnormal{(1)}] The series \eqref{omega} converges for $s \in \IC$ with $\mathrm{Re}(s) > \frac23$, and defines a holomorphic function in that region.
		\item[\textnormal{(2)}] The function $\omega$ can be analytically continued to a holomorphic function on $\IC \setminus \left( \{ \frac23 \} \cup (\frac12 - \N_0)\right)$. 
		\item[\textnormal{(3)}] The function $\omega$ has a simple pole at $s = \frac23$ with residue 
		\begin{align*}
		\mathrm{Res}_{s = \frac23} \omega(s) = \frac{\Gamma\left( \frac13 \right)^2}{2\sqrt{3}\pi}.
		\end{align*}
	 For each $m \in \N_0$, it has a simple pole at $s = \frac12 - m$ with residue
		\begin{align*}
			 \mathrm{Res}_{s=\frac12 - m}\omega(s) = \frac{(-1)^m}{16^m}  \binom{2m}{m} \zeta\left( \frac12 - 3m\right).
		\end{align*}
		\item[\textnormal{(4)}] Let $I \subset \IR$ be a compact interval. For all $\sigma \in I$, $\omega(\sigma + it)$ grows at most polynomially as $|t|\to\infty$, where the polynomial only depends on $I$.  
		\item[\textnormal{(5)}] We have  $\omega(-n) = 0$ for all $n \in \N$. 
	\end{enumerate}
\end{theorem}

We define for $\eta \in \IR \setminus \{ -\frac23, -\frac12, 0, \frac12, \frac32, \frac52, ...\}$ 
\begin{align*}
	E(\eta;z) := \frac{1}{2\pi i} \int_{-\eta-i\infty}^{-\eta+i\infty} J(s;z) ds,
\end{align*}
where
\[
	J(s;z) := \left(\frac2z\right)^s \Gamma(s)\zeta(s+1)\omega(s).
\]

The next proposition provides some basic facts about $E(\eta;z)$. As usual, $\mathrm{Log}$ denotes the prinicpal branch of the complex logarithm. For the proof, one uses Theorems \ref{GammaCollect}, \ref{ZetaFunc}, and \ref{omegapro} as well as results from \cite{Ro}.
 \begin{proposition} \label{ErrorHolo}
 	\ \begin{enumerate}[leftmargin=*]
 	\item[\textnormal{(1)}] Let $\eta \in \IR \setminus \{-\frac23, -\frac12, 0, \frac12, \frac32, ...\}$ and $0 < \delta \leq \frac{\pi}{2}$. As $z \to 0$ in $\mathcal{C}_\delta$, 
 	\begin{align*}
 	E(\eta;z) = O_{\eta, \delta}\left( |z|^\eta\right).
 	\end{align*}
 	\item[\textnormal{(2)}] Let $\eta \in \IR^+ \setminus \frac12(2\N_0+1)$. Then we have the functional equation 
 	\begin{equation*}
 	\label{Eequat} E(\eta;z) = E(\eta+1;z) + \mathrm{Res}_{s=-w} J(s;z),
 	\end{equation*}
 	where $w \in \N_0 + \frac 12$ is unique with the property $\eta < w < \eta + 1$. If $m+\frac12<\eta_1,\eta_2<m+\frac32$ for some $m \in \N_0$, then we have $E(\eta_1;z) = E(\eta_2;z)$ for all $z \in -i\mathbb{H}$. 
 	\item[\textnormal{(3)}] For all $\eta \in \IR \setminus \{-\frac23, -\frac12, 0, \frac12, \frac32, ...\}$, the function $z \mapsto E(\eta;z)$ is holomorphic on $-i\mathbb{H}$. For all  $\eta < -\frac23$ and $z \in -i\mathbb{H}$, we have 
 	\begin{equation*} \label{FullE}
 	\mathrm{Log}\left( G\left( e^{-z}\right)\right) = E(\eta; z).
 	\end{equation*}
 	\end{enumerate}
 	\end{proposition}

The asymptotic expansion of $\mathrm{Log}\left(G\left( e^{-z}\right)\right)$ can be shown similarly as in \cite{Ro}. 
\begin{theorem}\label{Asymptotic}
	Let $\eta>\frac12$ with $\eta\notin \frac12(2\N+1)$. Then for $z \in -i\mathbb{H}$ we have 
	\begin{multline*}
	\label{Heq}\mathrm{Log} \lp G\left( e^{-z}\right)\rp  
	= \frac{2^{\frac{2}{3}}3 X^{\frac{10}{3}}}{z^{\frac23}} - \frac{\sqrt{2}Y}{z^{\frac12}} - \frac13 \textnormal{Log}(z) + \frac 13 \log\left(16\pi^3\right) + z^{\frac 12} \sum_{0 \leq m < \eta - \frac12} \nu_m z^{m} 
	+ E\left(\eta;z\right),
	\end{multline*}
	where the coefficients $\nu_m$ have the explicit shape 
	\begin{equation*}\label{coeff}
	\nu_m = \frac{\sqrt{2\pi}}{(16\pi)^3 (8^5 \pi^4)^{m}} \frac{1}{m+1} \binom{2m}{m} \frac{(6m+6)!}{(3m+3)!} \zeta\left(m+\frac12\right) \zeta\left(3m+\frac72\right).
	\end{equation*}
\end{theorem}

\section{Error bounds}\label{Eb}

To give an asymptotic formula for $r(n)$, we  use Wright's Circle Method and the Saddle Point Method. For this we need many technical lemmas. In the following we always assume $\alpha \in \N + \frac14$. To find a satisfactory estimate for $E(\a;z)$ we define, with $k_n := \lfloor n^{\frac{2}{55}} \rfloor$, and $z \in -i\mathbb{H}$
\begin{align}
\label{E} \mathcal{E}_n(z) := E\left( k_n+\frac14;\frac{2X^2z}{n^{\frac35}}\right).
\end{align}

The goal of this section is to prove the following.

\begin{theorem}\label{Jsup}
	For $z \in \mathcal{C}_{\frac{\pi}{4}}$ we have, for some $B>0$
	\begin{align*}
	\mathcal{E}_n(z) \ll \left( \frac{B}{n^{\frac{1}{5}}}\right)^{n^{\frac{2}{55}}} |z|^{k_n + \frac14}.
	\end{align*}
\end{theorem}

To prove Theorem \ref{Jsup}, we need some auxiliary lemmas. The first proposition provides a uniform bound for $\zeta(s)$ on vertical lines. Using Theorem \ref{ZetaFunc} (1), Theorem \ref{GammaCollect} (1), (3), and Lemma \ref{prep} (1), we obtain.

\begin{proposition}\label{zetaest}
	For $\b \geq \frac14$ and $T\ge0$ we uniformly have
	\begin{align*}
	|\zeta(-\b + iT)| \ll (2\pi)^{-\b} \Gamma(\b + 1) \max\left\{1,T^{\b + \frac12}\right\}.
	\end{align*}
\end{proposition}

The following statement was proved in \cite{Ro}.

\begin{proposition}\label{Omegaexpr} Let $M \in \N$. We have for $s \in \IC$ with $\frac{3}{4} - \frac{M}{2} < \mathrm{Re}(s) < M+\frac12$
	\begin{multline*}
			\omega(s)  = \frac{ \Gamma(2s-1) \Gamma(1-s)\zeta(3s-1)}{\Gamma(s)} + \frac{1}{\Gamma(s)} \sum_{k=0}^{M-1} (-1)^k \frac{\Gamma(s+k)}{k!} \zeta(2s+k) \zeta(s-k) \\
	+ \frac{1}{2\pi i \Gamma(s)} \int_{M-\frac12 - i\infty}^{M-\frac12 +i\infty} \Gamma(s+z) \Gamma(-z) \zeta(2s+z) \zeta(s - z) dz.
	\end{multline*}
\end{proposition}

The next proposition gives a uniform estimate for the function $\omega$.

\begin{proposition}\label{Omegaestimate}
	For $\alpha \in \N + \frac14 $ and $T \geq 0$ we uniformly have, with some $C>0$ 
	\begin{align*}
	|\omega(-\alpha + iT)| \ll C^\alpha \Gamma(3\alpha + 3) \Gamma\left( 3\alpha + \frac{21}{4}\right)\max\left\{1, T^{5\alpha + \frac{11}{2}}\right\}.
	\end{align*} 
\end{proposition}

We split the technical proof of Proposition \ref{Omegaestimate} in several lemmas. The following lemma considers the first summand in Proposition \ref{Omegaexpr} and follows by a direct calculation using Theorem \ref{GammaCollect} (1), (3) and Proposition \ref{zetaest}.
\begin{lemma} \label{firstest} For $\alpha \in \N+\frac14$ and $T \geq 0$ we have uniformly
	\begin{align*}
	\left| \frac{\Gamma(-2\alpha-1+2iT) \Gamma(\alpha+1-iT)\zeta(-3\a-1 + 3iT)}{\Gamma(-\alpha + iT)} \right| \ll \left( \frac{9}{8\pi^3} \right)^\a \Gamma(3\alpha + 2) \max\left\{ 1, T^{3\alpha + \frac{3}{2}}\right\} e^{-\frac{\pi T}{2}}.
	\end{align*}
\end{lemma}

In the next lemma we give an estimate for the second summand in Proposition \ref{Omegaexpr}; the proof uses Proposition \ref{zetaest}, Theorem \ref{GammaCollect} (6), Lemma \ref{prep} (2), and Theorem \ref{ZetaFunc} (2).

\begin{lemma} \label{sumest} For $\alpha \in \N+\frac14$ and $T\ge0$, we uniformly have 
	\begin{align*}
	&\left| \sum_{k=0}^{2\lfloor\a\rfloor+ 2} (-1)^k \frac{\Gamma(-\a+k + iT)}{\Gamma(-\alpha + iT)k!}\zeta(-2\a+k+2iT) \zeta(-\a-k+ iT)\right| \\ 
	& \hspace{6cm} \ll \left(2\pi^3\right)^{-\a} \Gamma\left(3\alpha + \frac52\right) \max\left\{ 1, T^{5\alpha + \frac72}\right\}.
	\end{align*}
\end{lemma}

The following lemma deals with the integral in Proposition \ref{Omegaexpr}.

\begin{lemma} \label{integralest}
For $\a \in \N + \frac14$ and $T\geq0$ we uniformly have, for some $C>0$
\begin{align*}
& \left| \int_{2\alpha + 2 - i\infty}^{2\alpha + 2 +i\infty} \frac{\Gamma(-\alpha + iT+z) \Gamma(-z) \zeta(-2\a + 2iT+z) \zeta(-\a + iT - z)}{\Gamma(-\a + iT)} dz \right| \\
& \hspace{6cm} \ll C^\alpha \Gamma(3\alpha + 3) \Gamma\left( 3\alpha + \frac{21}{4}\right)\max\left\{1, T^{5\alpha + \frac{11}{2}}\right\}.
\end{align*} 
\end{lemma}
\begin{proof}
Choosing $M := 2\lfloor\a\rfloor + 3$ and $s := -\a + iT$, the term including integral in Proposition \ref{Omegaexpr} equals after substituting $z = 2\alpha + 2 + iv$
	\begin{align}
	\label{substitution} 	\frac{1}{2\pi}   \int_{- \infty}^{\infty} \frac{\Gamma\left(\alpha + 2+ i(T+v)\right) \Gamma\left(-2\alpha - 2 - iv\right) \zeta\left(2 + i(2T+v)\right) \zeta\left(-3\alpha - 2 + i(T-v)\right)}{\Gamma(-\alpha + iT)} dv. 
	\end{align}
 Now we have by Theorem \ref{GammaCollect} (1)
\begin{align}
\label{GammaUmformung} & \frac{\Gamma\left(\alpha + 2 + i(T+v)\right) \Gamma\left(-2\alpha - 2 - iv\right)}{\Gamma(-\alpha + iT)} \\ \nonumber
& \hspace{2cm} = \frac{\Gamma\left(\frac54 + i(T+v)\right) \Gamma\left(\frac12 - iv\right)}{\Gamma\left(\frac34+iT\right)} 
\frac{\prod_{j=0}^{\lfloor\a\rfloor}(-\alpha + j + iT) \prod_{j=0}^{\lfloor\a\rfloor} \left(\alpha + 1 - j + i(T+v)\right)}{\prod_{j=0}^{2\lfloor\a\rfloor + 2}\left(-2\alpha - 2 + j - iv\right)} .
\end{align}
We show by a direct calculation that uniformly for $v \in \IR$. We have, for some $B>0$,
\begin{align}
\label{vproduct} \left| \frac{\prod_{j=0}^{\lfloor\a\rfloor}(-\alpha + j + iT) \prod_{j=0}^{\lfloor\a\rfloor} \left(\alpha + 1 - j + i(T+v)\right)}{\prod_{j=0}^{2\lfloor\a\rfloor + 2}\left(-2\alpha - 2 + j - iv\right)} \right| \ll \frac{B^{2\alpha + \frac32} \max\left\{1,T^{2\alpha + \frac32}\right\}}{\max\{1, |v|\}}.
\end{align}

	To estimate the first factor in \eqref{GammaUmformung}, we distinguish two cases. First, assume $T + v \geq 0$. Then, with Theorem \ref{GammaCollect} \ref{Gamma}, \eqref{GammaUmformung}, and \eqref{vproduct} we obtain
	\begin{equation}\label{GammaFirst}
	 \left| \frac{\Gamma\left(\alpha + 2 + i(T+v)\right) \Gamma\left(-2\alpha - 2 - iv\right)}{\Gamma(-\alpha + iT)} \right| 
	\ll B^{2\alpha + \frac32} \max\left\{ 1,T^{2\alpha + 2 }\right\} \min\left\{ 1, e^{-\pi v}\right\}.
	\end{equation}

	If $T + v < 0$, then $v < -T \leq 0$, and hence we obtain with Theorem \ref{GammaCollect} (3), \eqref{vproduct}, and  \eqref{GammaUmformung}
	\begin{align}\label{GammaSecond}
	 \left| \frac{\Gamma\left(\alpha +2 + i(T+v)\right) \Gamma\left(-2\alpha - 2 - iv\right)}{\Gamma(-\alpha + iT)} \right|  \ll B^{2\alpha + \frac32}  \max\left\{ 1,T^{2\alpha + 2 }\right\} e^{\pi T} e^{\pi v}.
	\end{align}
	We conclude with $|\zeta(2 + i(2T+v))| \leq \zeta(2)$, using \eqref{GammaFirst} and \eqref{GammaSecond} that \eqref{substitution} can be bound against a constant times 
	\begin{align}\label{greatestimate}
		 B^{2\alpha + \frac32} \max\left\{ 1,T^{2\alpha + 2 }\right\} \left(e^{\pi T} I_1(T) + I_2(T)\right),
	\end{align}
where 
\begin{align*}
I_1(T)&:=\int_{-\infty}^{-T}  \left| \zeta\left(-3\alpha - 2 + i(T-v)\right) \right| e^{\pi v} dv ,
\\
 I_2(T)&:=\int_{-T}^\infty  \left| \zeta\left(-3\alpha - 2 + i(T-v)\right) \right| \min\left\{ 1, e^{-\pi v}\right\} dv.
\end{align*}
 We consider  $I_1(T)$ and $I_2(T)$ separately. We obtain substituting $u = T-v$ into $I_1(T)$
	\begin{align*}
	& e^{\pi T} I_1(T) = 	e^{2\pi T}  \int_{2T}^{\infty} \left| \zeta\left(-3\alpha - 2 + iu\right) \right| e^{-\pi u} du.
\end{align*}
	With Proposition \ref{zetaest} we can  bound  
			\begin{align} \label{zwischen}
	e^{\pi T} I_1(T) 
	 \ll \left(2\pi^2\right)^{-3\alpha} \Gamma\left( 3\alpha +3\right) \Gamma\left(3\alpha + \frac72; 2\pi T\right) e^{2\pi T} + (2\pi)^{-3\alpha} \Gamma\left( 3\alpha + 3 \right).
	\end{align}
Using the monotonicity of the $\Gamma$-function to bound
	$
	\Gamma(3\alpha + \frac72;2\pi T) \ll \Gamma(3\alpha + \frac{17}{4};2\pi T),
	$ (so that $3\a+\frac{17}{4}\in\N$)
	 we use Theorem \ref{IncompleteGamma} (2) to give 
	\begin{align*}
	\Gamma\left(3\alpha + \frac{17}{4};2\pi T\right) 
	\ll (2\pi)^{3\alpha} \Gamma\left(3\alpha + \frac{21}{4}\right)  \max\left\{1, T^{3\alpha + \frac{13}{4}}\right\} e^{-2\pi T}.
	\end{align*}
Plugging this in \eqref{zwischen} we conclude 
	\begin{align} 
	\label{firstInt}  B^{2\alpha + \frac32} \max\left\{ 1,T^{2\alpha + 2 }\right\} e^{\pi T} I_1(T) 
	\ll B^{2\alpha + \frac32} \Gamma\left(3\alpha + \frac{21}{4}\right) \Gamma\left( 3\alpha + 3\right)\max\left\{1, T^{5\alpha + \frac{21}{4}}\right\}.
	\end{align}
	
	In $I_2(T)$ we split the integral at $v=0$.
	First, with Proposition \ref{zetaest} we obtain, 
	\begin{align}	\label{I2first}
	 \int_{-T}^0  \left| \zeta\left(-3\alpha - 2 + i(T-v)\right) \right| dv \ll \pi^{-3\alpha}\Gamma\left( 3\alpha +3\right) \max\left\{ 1, T^{3\alpha +\frac72}\right\}.
	\end{align}
	
	For the contribution from $v\ge0$, we substitute $u = T-v$. For the contribution from $u\le0$, we obtain with Proposition \ref{zetaest} 
	\begin{align}\label{I2sec1}
	 e^{-\pi T} \int_{-\infty}^{0} \left| \zeta\left(-3\alpha - 2 + iu\right) \right| e^{\pi u} du \ll \left(2\pi^2\right)^{-3\alpha} \Gamma\left( 3\alpha + 3 \right)  \Gamma\left(3\alpha + \frac72\right) e^{-\pi T}.
\end{align}
For the contribution from $0\le u\le T$ we obtain again with Proposition \ref{zetaest}
\begin{align}
 \label{I2sec2} e^{-\pi T} \int_{0}^{T} \left| \zeta\left(-3\alpha - 2 + iu\right) \right| e^{\pi u} du \ll (2\pi)^{-3\alpha} \Gamma\left( 3\alpha + 3 \right)  \max\left\{1, T^{3\alpha + \frac52}\right\}.
\end{align}
As a result, we obtain with \eqref{I2first}, \eqref{I2sec1}, and \eqref{I2sec2} the estimate 
\begin{align*}
 I_2(T)  \ll \Gamma\left( 3\alpha + 3 \right)  \Gamma\left(3\alpha + \frac72\right) \max\left\{ 1, T^{3\alpha +\frac72}\right\}.
\end{align*} 
	Hence we obtain 
	\begin{align}\label{secInt}
	   B^{2\alpha + \frac32} \max \left\{ 1,T^{2\alpha + 2 } \right\} I_2(T) 
	  \ll B^{2\alpha + \frac32} \Gamma\left( 3\alpha + 3 \right)  \Gamma\left(3\alpha + \frac72\right) \max\left\{ 1, T^{5\alpha +\frac{11}{2}}\right\}.
	\end{align}
Together with \eqref{greatestimate}, \eqref{firstInt}, and \eqref{secInt} we conclude the claim.
\end{proof}
Now, we are ready to prove Proposition \ref{Omegaestimate}.

\begin{proof}[Proof of Proposition \ref{Omegaestimate}]
	Using Proposition \ref{Omegaexpr}, Lemmas \ref{firstest}, \ref{sumest}, and \ref{integralest} the claim follows.
\end{proof}

We are now ready to prove Theorem \ref{Jsup}.

\begin{proof}[Proof of Theorem \ref{Jsup}] Using the triangle inequality we obtain 
\[
|\mathcal{E}_n(z)| \ll \int_{-\infty}^\infty \left| \left(  \frac{n^{\frac35}}{X^2 z} \right)^{-k_n - \frac14 + iv} \Gamma\left(-k_n - \frac14 + iv\right) \zeta\left( -k_n + \frac34 + iv \right)  \omega\left( -k_n - \frac14 + iv \right)  \right| dv.
\]
With Theorem \ref{GammaCollect} (4), Theorem \ref{ZetaFunc} (1), and Proposition \ref{Omegaestimate}, we obtain
	\begin{multline*}
	\Gamma\left(-k_n - \frac14 + iv\right) \zeta\left( \frac{3}{4} - k_n + iv\right) \omega\left(-k_n - \frac14 + iv\right) \\ 
	\ll C^{k_n} \Gamma\left(3k_n + \frac{15}{4}\right) \Gamma\left( 3k_n+ 6\right)\max\left\{1, |v|^{5k_n + 7}\right\}e^{-\frac{\pi |v|}{2}}
	\end{multline*}
	 for some constant $C > 0$. It follows for some constant $C_1 > 0$ 
	\begin{equation}\label{Eest}
		|\mathcal{E}_n(z)| \ll C_1^{k_n} n^{-\frac{3}{5}\left(k_n + \frac14\right)} |z|^{k_n+\frac14}  \Gamma\left(3k_n + 4\right) \Gamma\left( 3k_n+ 6\right) \Gamma\left( 5k_n + 8\right).
	\end{equation}
Since
	\begin{equation*}
	 	\Gamma(3k_n + 4) \Gamma(3k_n+ 6) \Gamma(5k_n + 8)  \ll k_n^{15} \Gamma(3k_n+1)^2\Gamma(5k_n + 1),
	\end{equation*}
	we obtain with Theorem \ref{GammaCollect} (2) for constants $C_2, C_3 > 0$
	\begin{equation}\label{GammaStir}
		\Gamma(3k_n + 4) \Gamma(3k_n+ 6) \Gamma(5k_n + 8) \ll n^{C_2} C_3^{k_n} n^{\frac{2}{5}n^{\frac{2}{55}}}.
	\end{equation}
	Note that $n^{-\frac{3k_n}{5}} \ll n^{\frac{3}{5}} n^{-\frac{3}{5}n^{\frac{2}{55}}}$ since $k_n = \lfloor n^{\frac{2}{55}} \rfloor \geq n^{\frac{2}{55}} - 1$. Hence, with \eqref{Eest} and \eqref{GammaStir},
\begin{align}
\label{GammaStir2} n^{-\frac{3}{5}\left(k_n + \frac14\right)} \Gamma(3k_n + 4) \Gamma(3k_n+ 6) \Gamma(5k_n + 8) \ll n^{C_2+\frac{9}{20}} \left( \frac{C_3}{n^{\frac{1}{5}}}\right)^{n^{\frac{2}{55}}}. 
\end{align}
Since $n^{C_2+\frac{9}{20}} \ll C_4^{n^{\frac{2}{55}}}$ for some $C_4 > 1$,  we finally obtain with $B := C_1C_3C_4$, $B^{k_n} = B^{\lfloor n^{\frac{2}{55}}\rfloor} \ll B^{n^{\frac{2}{55}}}$, \eqref{Eest}, and \eqref{GammaStir2} the theorem. 
\end{proof}

\section{The saddle point function}
 
 When using the Saddle Point Method we approximate the major arc integral \eqref{MajorFinal} below to obtain an asymptotic expansion of $r(n)$ as $n \rightarrow \infty$. The major arc integrals have the form
\begin{align}
\label{abstractsaddle} I_n := \int_{1-i}^{1+i} G_n(z)e^{n^{\frac{2}{5}} F_n(z)} dz.
\end{align}
 The $G_n$ are holomorphic functions that are uniformly bounded on compact subsets of the right half-plane (see Lemmas \ref{hsup} and \ref{Jestimate}) and 
\begin{align*} 
F_n(z) := \frac{3X^2}{z^{\frac{2}{3}}} - \frac{Y}{Xn^{\frac{1}{10}} z^{\frac12}} + 2X^2z.
\end{align*}
To use the Saddle Point Method one needs to find saddle points which play the key role to find an asymptotic expansion for \eqref{abstractsaddle}; which then is achieved in Lemma \ref{SaddlePoint} below. In the following we sketch the ideas for the proof. We recall that a saddle point $z_0$ satisfies $F'_n(z_0) = 0$. However, since the saddle point depends on $n$ it is therefore useful to think instead of a ``saddle point function'' (explained below). For this consider the function 
\begin{align*}
 \mathcal{F}(z; w) := -\frac{2X^2}{z^{\frac{5}{3}}} + \frac{Yw}{2X z^{\frac{3}{2}}} + 2X^2,
\end{align*}
which is holomorphic on $-i\mathbb{H} \times \IC$; (note that $F'_n(z) = \mathcal{F}(z; n^{-\frac{1}{10}})$). When $w$ is fixed and small, we find below that $\mathcal{F}$ has an inverse $z \mapsto \mathcal{F}^{-1}(z; w)$ (i.e., with respect to the variable $z$) locally around $\mathcal{F}(1;w) = \frac{Yw}{2X}$. It is again holomorphic and possesses a power series expansion $\mathcal{F}^{-1}(z; w)=1 + \sum_{j=1}^\infty a_w(j)(z - \frac{Yw}{2X})^j$. We would like to obtain roots of the equations $F_n'(z_0) = 0$ by just plugging in $z = 0$ and $w= n^{-\frac{1}{10}}$ in this power series, but it is not immediately clear that the radius of convergence is large enough to guarantee convergence. We remedy this by adding another component in the image and looking at the function 
\begin{equation*}
\mathcal{H}(z;w) := (\mathcal{F}(z;w), w) .
\end{equation*}
To this function we can then apply the Inverse Mapping Theorem to give a local holomorphic inverse $\mathcal{H}^{-1}$ around $(0,0)$ and consider pre-images $\mathcal{H}^{-1}(0;x) =: (S(x),x)$ for $x$ real and small. This defines a real-valued function $x \mapsto S(x)$ in a small neighborhood of $0$, that we call {\it saddle point function}. We then have
\begin{align*}
 \mathcal{F}(S(x);x) = 0.
\end{align*}
The most important case for us, $x = n^{-\frac{1}{10}}$, is denoted by 
\begin{align}
\label{Sn} S_n := S\left(n^{-\frac{1}{10}}\right). 
\end{align}
In particular, $F_n'(S_n) = 0$ for $n$ sufficiently large. In the next lemma we prove properties of $S(x)$.

\begin{lemma} \label{SaddlePoint} The function $S$ is holomorphic near $0$, and has a power series expansion
	\begin{align}\nonumber
		S(x) &= 1 + \sum_{m=1}^\infty \varrho(m)x^m\\
		\label{SDef}
		&= 1 - \frac{3Y}{20X^3}x - \frac{3Y^2}{800X^6}x^2 - \frac{11Y^3}{64000X^9}x^3 + \frac{4959Y^5}{2048000000X^{15}} x^5 +O\left(x^6\right),
	\end{align}
	with some $\varrho(m)\in\IR$. In particular, $x \mapsto S(x)$, for $x$ sufficiently small, is a real-valued function and the sequence $S_n$ converges to $1$ and satisfies $S_n < 1$ for $n$ sufficiently large. 
\end{lemma}

\begin{proof}
	Note that $\mathcal{H}$ is holomorphic in $-i\mathbb{H} \times \IC$. We obtain for the Jacobi matrix
	\begin{equation*}\label{Jacoby}
	J_{\mathcal{H}}(z;w) := \mat{\frac{\partial}{\partial z} \mathcal{F}(z; w) & \frac{\partial}{\partial w} \mathcal{F}(z; w) \\ \vspace{-.15cm}&\vspace{-.15cm} \\ \frac{\partial}{\partial z} w & \frac{\partial}{\partial w} w} = \mat{ \frac{10 X^2}{3 z^{\frac{8}{3}}}-\frac{3 w Y}{4 X z^{\frac{5}{2}}} & \frac{Y}{2X z^{\frac{3}{2}}} \\ \vspace{-.15cm}&\vspace{-.15cm} \\ 0 & 1 }.
	\end{equation*}
	In particular $\det(J_{\mathcal{H}}(1;0))=\frac{10}{3}X^2 \not= 0$. Therefore, due to the Inverse Mapping Theorem (see \cite{Kaup}, p. 27), $\mathcal{H}$ is locally biholomorphic around $(1,0)$ with inverse $\mathcal{H}^{-1}=:(\mathcal{H}_1^{-1},\mathcal{H}_2^{-1})$, i.e., $\mathcal{H}^{-1}$ can be expanded locally into a power series around $\mathcal{H}(1;0) = (0,0)$. Hence, we obtain 
	\begin{align*}
	\mathcal{H}^{-1}_1(u;w) = \sum_{k, m \ge 0} \beta(k,m) u^k w^m, \qquad \mathcal{H}^{-1}_2(u;w) = w.
	\end{align*} 
	As a result, for real $x$ sufficiently small, we find $\mathcal{H}^{-1}(0;x)=(S(x),x)$. Hence we obtain 
	\begin{align*}
	S(x) = \sum_{m=0}^\infty \beta(0,m)x^m = \sum_{m=0}^\infty \varrho(m)x^m,
	\end{align*}
	and $\varrho(0) = 1$, since we have $\mathcal{H}^{-1}(0;0) = (1,0)$. Next note that all $\varrho(m)$ are real. Indeed this holds, since $\overline{\mathcal{F}(z;x)} = \mathcal{F}(\overline{z};x)$, both $S(x)$ and $\overline{S(x)}$ are zeros of $z \mapsto \mathcal{F}(z;x)$; because of biholomorphicity of $\mathcal{H}$ we obtain $S(x) = \overline{S(x)}$. One can use a computer to find the first few coefficients in \eqref{SDef}. A straightforward calculation using \eqref{Sn} and \eqref{SDef} shows that $S_n \to 1$, as $n \to \infty$, and $S_n < 1$ for all $n$ sufficiently large. 
	\end{proof}
  
\section{Approximation of holomorphic functions}\label{sec:appholf}

\subsection{Error bounds for the asymptotic terms of $G$}
In the next proposition we give a criterion for monotonicity of power series; we leave the straightforward proof to the reader.

\begin{proposition} \label{PowerHelp} Let $a(n)$ be a decreasing sequence of non-negative real numbers that converges monotonously to 0, such that $a(0) \geq a(1) > a(2) \geq 0$. Then the function $P : (-1, 1) \to \IR$.
	\begin{align*}
	P(x) := \sum_{n=0}^\infty (-1)^n a(n) x^{2n}
	\end{align*}
	has the global maximum $a(0)$ at $x=0$, is increasing in $[-2^{-\frac12},0]$, and decreasing in $[0,2^{-\frac12}]$.
\end{proposition}

In the next lemma we investigate the following function that is important in the Saddle Point Method. We define, with $S_n$ the sequence of saddle points provided in Lemma \ref{SaddlePoint},
\begin{align}
\label{Fn} f_n(x) := \frac{3X^2}{(S_n + ix)^{\frac23}}  - \frac{Y}{Xn^{\frac{1}{10}}(S_n + ix)^{\frac12}} + 2X^2(S_n+ix).
\end{align}
\begin{lemma} \label{Mono} For all $n$ sufficiently large, $x \mapsto \mathrm{Re}(f_n(x))$ has its global maximum at $x=0$. In $[-\frac12, 0]$ it is monotonously increasing, in $[0, \frac12]$ it is monotonously decreasing, and for $|x| \geq \frac12$ and real $\kappa$ sufficiently small   
	\begin{align}
	\label{rEsti} |\mathrm{Re}(f_n(x + i\kappa))| \leq 4.8X^2.
	\end{align}
\end{lemma}

\begin{proof} 
	Without loss of generality we can assume $x \geq 0$. Assume that $n$ is sufficiently large such that $\frac12 < S_n < 1$. Let $x \in (0, \frac12]$. We compute 
	\begin{align}
	\label{rFormula} \mathrm{Re}(f_n(x)) = \frac{3X^2}{S_n^{\frac{2}{3}}} \sum_{j=0}^\infty \binom{- \frac23}{2j} (-1)^j \left( \frac{x}{S_n}\right)^{2j} - \frac{Y}{Xn^{\frac{1}{10}}S_n^{\frac12}} \sum_{j=0}^\infty \binom{-\frac12}{2j} (-1)^j \left( \frac{x}{s_n}\right)^{2j} + 2X^2S_n.
	\end{align}
	Note that for $n$ sufficiently large, one can show, using Lemma \ref{prep} (3), that the sequence 
	\begin{align*}
	u_n(j) := \frac{3X^2}{S_n^{\frac{2}{3}}} \binom{- \frac23}{2j} - \frac{Y}{Xn^{\frac{1}{10}}S_n^{\frac12}} \binom{-\frac12}{2j}   
	\end{align*}
	satisfies $u_n(j) > u_n(j+1) > 0$ for all $j \in \N$ and $u_n(0) > u_n(1)$. With \eqref{rFormula}, $u_n(0) > u_n(1) > u_n(2) > 0$, and Lemma \ref{PowerHelp} one concludes monotonicity and maximality of $\mathrm{Re}(f_n)$ at $x=0$. 
	
	The estimate \eqref{rEsti} follows by a straightforward calculation.
\end{proof}

Using Lemma \ref{Mono} we can show that the exponentials in the major arc integral \eqref{MajorFinal} are negligible on a specific part of the integration curve. 
\begin{lemma}\label{Saddlesup}
	We have 
	\begin{multline*}
	\hspace{-.25cm}\sup_{- S_n \leq x \leq -n^{- \frac{7}{40}}} \left| \exp\left( n^{\frac{2}{5}} f_n\left( x \right) \right) \right| \ll \exp\left( n^{\frac{2}{5}} \left( \frac{3X^2}{S_n^{\frac23}} + 2X^2 S_n - \frac{Y}{Xn^{\frac{1}{10}} S_n^{\frac12}} - \frac{5X^2}{3n^{\frac{7}{20}}S_n^{\frac{8}{3}} } \right) \right), \quad n \rightarrow \infty.
	\end{multline*}
\end{lemma}

The following key lemma is due to Romik and helps us to obtain estimates on the minor arcs. 

\begin{lemma}[Romik \cite{Ro}, equation (94)] \label{Gestimate} For all $\kappa > 0$ there exist constants $\beta, \delta > 0$, such that for all $0 < t < \beta$ and $\kappa t \leq |u| \leq \pi$ we have 
	\begin{align*}
	\left|G\left(e^{-t+iu}\right)\right| \leq \exp\left( - \frac{\delta}{\sqrt{t}}\right) G\left(e^{-t}\right).
	\end{align*}
\end{lemma}

\subsection{Approximation for the auxillary asymptotic terms in $G$}
 Using Theorem \ref{GammaCollect} (2)
 we obtain by a straightforward calculation the following asymptotic bound for the coefficients $\nu_m$ defined in  Theorem \ref{Asymptotic}.
\begin{lemma}\label{nuestimate}
	There exists a constant $C > 0$, such that we have for $m \in \N$ 
	\begin{align*}
	\nu_m \ll C^m m^{3m}.
	\end{align*}
\end{lemma}

To find the asymptotic behavior of $r(n)$, by Theorem \ref{Asymptotic} we investigate the behavior of $G(e^{-z})$ towards $q=1$, where we choose $k_n - \frac12 < \eta_n < k_n + \frac12$  in Theorem \ref{Asymptotic}. In Section 4 the error integral from Theorem \ref{Asymptotic} is already bounded. The main term, i.e., the terms with negative exponents of the expansion in Theorem \ref{Asymptotic} is investigated in Lemmas \ref{Mono} and \ref{Saddlesup}. We now look at the remaining terms with positive exponents. Consider the following function $H_n : B_{\frac12}(0) \times B_\delta(0) \to \IC$
\begin{align}
\label{Hn} H_n(w;z) := \exp\left(\sum_{m=0}^{k_n-1} \nu_m \left( 2X^2 (S(z) + iw) \right)^{m+\frac12} z^{6m+3}\right),
\end{align}
where $S(z)$ is the saddle point function \eqref{SDef}, that is holomorphic in some region $B_\delta(0)$. When choosing $\delta$ sufficiently small, we achieve $S(B_\delta(0)) \subset B_{\frac14}(1)$, since $S(0) = 1$. Hence, for all $z \in B_\delta(0)$ and $w \in B_{\frac12}(0)$, we obtain $\mathrm{Re}(S(z) + iw) > \frac{3}{4} - \frac12 = \frac14$. Since the principal branch of the square root is holomorphic on the right half-plane, $H_n$ is holomorphic on $B_{\frac12}(0) \times B_\delta(0)$ as a composition of holomorphic functions. As such, we can write $H_n(w;z)$ as a power series 
\begin{align}
\label{cn} H_n\left(w; z\right) = \sum_{m=0}^\infty c_{n,z}\left(m\right) w^m = \sum_{m,\ell \geq 0} a_{n,m}(\ell)z^\ell w^m. 
\end{align}
Our goal is to show the following proposition.  
\begin{proposition}\label{HFinal}
	Fix $M\in\N_0$. For $|x| \leq n^{\frac{1}{40}}$ and $n$ sufficiently large we have
	\begin{align*}
		H_n\left( \frac{x}{n^{\frac{1}{5}}}; n^{-\frac{1}{10}}\right) = 1 + \sum_{m=1}^M \frac{P^{[1]}_m(x)}{n^{\frac{m}{10}}} + O_M\left( n^{-\frac{M+1}{15}}\right),
	\end{align*}
	where $P^{[1]}_m$ are polynomials of degree at most $\frac{m}{2}$ that do not depend on $n$. 
\end{proposition}

To prove Proposition \ref{HFinal}, we need the following lemma. 

\begin{lemma} \label{hsup}  Let $0<\delta < 1$ be fixed. The sequence $w\mapsto H_n( w; n^{-\frac{1}{10}})$ of holomorphic functions converges compact to 1 on $\mathbb{C}_\delta$.
\end{lemma}

\begin{proof}
	 Since by Lemma \ref{SaddlePoint} the sequence $S_n=S(n^{-\frac{1}{10}})$ converges to 1 as $n \rightarrow \infty$, we have for all $w \in \mathbb{C}_\delta$, $\mathrm{Re}(S_n + iw) > 0$ for $n$ sufficiently large. For these $n$ the functions $w \mapsto H_n( w; n^{-\frac{1}{10}})$ are holomorphic in $\mathbb{C}_\delta$. Let $D \subset \mathbb{C}_\delta$ be compact. Then there exists a constant $C_D > 0$, such that we have $|S_n + iw| \leq S_n + |w| \leq C_D$ for all $z \in D$ and for all $n$ sufficiently large.  Using this, we obtain with the triangle inequality 
	\begin{align*}
		\left| H_n\left( w; n^{-\frac{1}{10}}\right) \right| & \leq \exp\left( \sum_{m=0}^{k_n - 1} \left| \nu_m \left( \frac{2C_DX^2}{n^{\frac{3}{5}}}  \right) \right|^{m+\frac12} \right).
	\end{align*}
	The claim now follows by Lemma \ref{nuestimate} completing the geometric series. 
\end{proof}

The next step is to study the coefficients $a_{n,m}(\ell)$ in \eqref{cn}. 

\begin{lemma} \label{independence} For $\ell \in \N_0$, there exists $N_{\ell} \in \N$, such that for all $m \in \N_0$, $a_{n,m}(\ell)$ is constant for $n \ge N_{\ell}$. 
\end{lemma} 
\begin{proof} 
The sequence $k_n = \lfloor n^{\frac{2}{55}} \rfloor$ is mononotous, and as a result we see for $n_2 \geq n_1$ 
	\begin{align}
	\label{Hdiff} H_{n_2}(w,z) - H_{n_1}(w,z) = O_{n_1, n_2, w}\left(z^{6k_{n_1}+3}\right).
	\end{align}
	On the other hand, we obtain with \eqref{cn}
	\begin{align}
	\label{HdiffPower} H_{n_2}(w,z) - H_{n_1}(w,z) = \sum_{m=0}^\infty \sum_{\ell=0}^\infty \left( a_{n_2,m}(\ell) - a_{n_1, m}(\ell)\right)z^\ell w^m.
	\end{align}
	From \eqref{Hdiff} and \eqref{HdiffPower} it follows, that $a_{n_2,m}(\ell) - a_{n_1, m}(\ell) = 0$ for all $m \in \N_0$ if $\ell < 6k_{n_1} + 3$. With $N_{\ell} := \ell^{55}$ the lemma follows. 
\end{proof}
Next, it is useful to have an uniform upper bound for the values $a_{n,m}(\ell)$. Using the Cauchy Integral formula we next obtain the following bound. 
\begin{lemma} \label{abound} We have for $m, \ell \in \N_0$ and $n$ sufficiently large
	\begin{align*}
	|a_{n,m}(\ell)| \ll 4^m n^{\frac{\ell}{30}}.
	\end{align*}
\end{lemma}

We next approximate the coefficients $c_{n,z}$, defined in \eqref{cn}.
\begin{lemma} \label{clarge} Fix $M\in\N_0$. Then for $n$ sufficiently large we have 
	\begin{align*}
	c_{n,n^{-\frac{1}{10}}}\left(m\right) = \sum_{\ell=0}^M \frac{a_m(\ell)}{n^{\frac{\ell}{10}}} + O\left( \frac{4^m}{n^{\frac{M+1}{15}}}\right),
	\end{align*}
	where the $a_m(\ell)$ are independent on $n$, and the $O$-constant is independent of $m$, $M$, and $n$.
	\end{lemma}
\begin{proof}
	By Lemma \ref{independence}, for each $0 \leq \ell \leq M$ there exists a constant $N_\ell$ such that $a_{n,m}(\ell)$ is constant for $n \geq N_{\ell}$. Since $M$ is fixed and independent on $n$, $R_{M} := \max_{0 \leq \ell \leq M} N_\ell$ is also independent on $n$. 
	Thus $a_{n,m}(\ell)$ does not depend on $n$ if $n \geq R_M$, and for these $n$ can simply be noted as $a_{m}(\ell)$. Hence, for all $n \geq R_M$ we have, by \eqref{cn} 
\begin{align*}
c_{n,z}(m) = \sum_{\ell=0}^M a_{m}(\ell) z^\ell + \sum_{\ell=M+1}^\infty a_{n,m}(\ell) z^\ell.
\end{align*}
	With Lemma \ref{abound} and \eqref{cn} we then obtain for $n$ sufficiently large 
	\begin{align*}
	\left| \sum_{\ell=M+1}^\infty \frac{a_{n,m}(\ell)}{n^{\frac{\ell}{10}}}\right| \ll 4^m \sum_{\ell=M+1}^\infty 
 n^{-\frac{\ell}{15}}
		 = 4^m \frac{n^{-\frac{M+1}{15}}}{1 - n^{-\frac{1}{15}}} = O\left( \frac{4^m}{n^{\frac{M+1}{15}}}\right).
	\end{align*}
	This proves the lemma. 
	\end{proof}
By \eqref{cn}, the Cauchy Integral Formula, Lemmas \ref{clarge} and \ref{hsup} Proposition \ref{HFinal} follows. 

\subsection{Approximation for the main asymptotic term in $G$}
The next step is to investigate the main term in Theorem \ref{Asymptotic}. Similarly as in the case of $H_n$, we first interpret the saddle point parameter as an independent variable of a holomorphic function. Therefore, we put for $w \in \mathbb{C}_{\frac12}$ and $z \in B_\delta(0)$, where $\delta$ is sufficiently small, 
	\begin{align*}
f(w,z) := \exp\left( z^{-4}\left( \frac{3X^2}{(S(z) + iw)^{\frac{2}{3}}} - \frac{Yz}{X(S(z) + iw)^{\frac12}} + 2X^2 (S(z) + iw)\right)\right).
\end{align*}
Moreover define
\begin{align*}
	A(n) := \exp\left( A_1n^{\frac{2}{5}} - A_2n^{\frac{3}{10}} - A_3n^{\frac{1}{5}} - A_4n^{\frac{1}{10}} \right).
\end{align*}

\begin{lemma}\label{fexpand}
	For $z$ and $x$ sufficiently small we have the Laurent expansion 
	\begin{multline}\label{separatez4}
		z^{-4}\left( \frac{3X^2}{(S(z) + ixz^2)^{\frac{2}{3}}} - \frac{Yz}{X(S(z) + ixz^2)^{\frac12}} + 2X^2 \left(S(z) + ixz^2\right)\right) \\
		= \frac{A_1}{z^4} - \frac{A_2}{z^3} - \frac{A_3}{z^2} - \frac{A_4}{z} - A_5 - \frac{5X^2}{3}x^2 + \sum_{m=1}^\infty P_{\ell+4}(x)z^\ell,
	\end{multline}
	where the $P_{\ell+4}$ are polynomials. For every fixed $M\in\N_{>1}$ and $|x| \leq n^{\frac{1}{40}}$, we have the expansion 
	\begin{align*}
		\exp\left( n^{\frac{2}{5}} f_n\left( \frac{x}{n^{\frac{1}{5}}}\right) \right) = A(n)\exp\left(- A_5\right) \exp\left( - \frac{5X^2}{3}x^2 \right) \left( 1 + \sum_{m=1}^M \frac{P^{[2]}_m(x)}{n^{\frac{m}{10}}} + O\left( n^{-\frac{3(M+1)}{80}} \right)\right),
	\end{align*}
	where $n$ is sufficiently large and $P^{[2]}_m(x)$ are polynomials independent from $n$ with $\deg(P^{[2]}_m) \leq 2m$ . 
	\end{lemma}
\begin{proof} 
	 We start with the Laurent expansion 
	\begin{align*}
	 z^{-4}\left( \frac{3X^2}{(S(z) + iw)^{\frac{2}{3}}} - \frac{Yz}{X(S(z) + iw)^{\frac12}} + 2X^2 \left(S(z) + iw\right)\right) =: z^{-4} \sum_{\ell, m \ge 0} a(\ell, m)w^\ell z^m .
	\end{align*}
	Setting $w = xz^2$ gives
	\begin{multline*} 
	z^{-4}\left( \frac{3X^2}{(S(z) + ixz^2)^{\frac{2}{3}}} - \frac{Yz}{X(S(z) + ixz^2)^{\frac12}} + 2X^2 \left(S(z) + ixz^2\right)\right) \\
	= z^{-4} \sum_{\ell,m\ge0} a(\ell, m)x^\ell z^{2\ell} z^m = z^{-4} \sum_{\ell=0}^\infty P_\ell(x)z^\ell,
	\end{multline*}
	where the polynomials $P_\ell$ are independent from $x$ and $z$ and satisfy $\deg(P_\ell) \leq \lfloor \frac{\ell}{2}\rfloor$. We obtain with the help of a computer
	\begin{align*}
	 \sum_{\ell=0}^\infty P_\ell(x)z^{\ell-4} = \frac{A_1}{z^4} - \frac{A_2}{z^3} - \frac{A_3}{z^2} - \frac{A_4}{z} - A_5 - \frac{5X^2}{3}x^2 + \sum_{\ell=5}^\infty P_\ell(x)z^{\ell-4}.
	\end{align*} 
	which proves \eqref{separatez4}. As a result, we obtain 
	\begin{align*}
	& f\left(\frac{x}{n^{\frac15}},n^{-\frac{1}{10}}\right) = A(n)\exp\left(- A_5 - \frac{5X^2}{3}x^2 + \sum_{\ell=1}^\infty \frac{P_{\ell+4}(x)}{n^{\frac{\ell}{10}}}\right).
	\end{align*}
	We write
	\begin{align*}
	 \exp\left( \sum_{\ell=1}^\infty \frac{P_{\ell+4}(x)}{n^{\frac{\ell}{10}}} \right) 
	=: 1 + \sum_{m=1}^M \frac{P^{[2]}_m(x)}{n^{\frac{m}{10}}} + \sum_{m=M+1}^\infty \frac{P^{[2]}_m(x)}{n^{\frac{m}{10}}},
	\end{align*}
	and we are left to show that uniformly for $|x| \leq n^{\frac{1}{40}}$
	\begin{align*}
		 \sum_{m=M+1}^\infty \frac{P^{[2]}_m(x)}{n^{\frac{m}{10}}} = O\left( n^{-\frac{3(M+1)}{80}}\right),
	\end{align*}
and that $\deg(P^{[2]}_m(x))\leq 2m$. This follows by a lengthy however straightforward calculation.
\end{proof}

In a similar way as demonstrated in Lemma \ref{fexpand} we show the following. 

\begin{lemma} \label{cuberootexpand} Fix $M\in \N_{>1}$. For all $n$ sufficiently large and $|x| \leq n^{\frac{1}{40}}$ we have
	\begin{align*}
	\left( S_n + \frac{ix}{n^{\frac{1}{5}}}\right)^{-\frac13} = 1 + \sum_{m=1}^M \frac{P^{[3]}_m(x)}{n^{\frac{m}{10}}} + O\left( n^{-\frac{7(M+1)}{80}}\right),
	\end{align*}
	where the $P^{[3]}_m$ are polynomials of degree at most $\frac{m}{2}$. 
\end{lemma}
\subsection{Polynomial approximations}
Now define (considered as formal series)
\begin{align*}
	 \prod_{j=1}^3 \left( 1 + \sum_{m=1}^\infty \frac{P^{[j]}_m(x)}{n^{\frac{m}{10}}} \right) =: \sum_{m=0}^\infty \frac{P^{[4]}_m(x)}{n^{\frac{m}{10}}}.
\end{align*}
We turn this into an analytic formula by truncating the series on the left-hand side ($M\in\N$)
\begin{align*}
\prod_{j=1}^3 \left( 1 + \sum_{m=1}^M \frac{P^{[j]}_m(x)}{n^{\frac{m}{10}}} \right) =: \sum_{m=0}^M \frac{P^{[4]}_m(x)}{n^{\frac{m}{10}}} + \sum_{m=M+1}^{3M} \frac{P^{[4]}_{M,m}(x)}{n^{\frac{m}{10}}}.
\end{align*}
Note that for $k \in \N$
\begin{align*}
\prod_{j=1}^3 \left( 1 + \sum_{m=1}^{M+k} \frac{P^{[j]}_m(x)}{n^{\frac{m}{10}}} \right) - \prod_{j=1}^3 \left( 1 + \sum_{m=1}^{M} \frac{P^{[j]}_m(x)}{n^{\frac{m}{10}}} \right) = O\left(n^{-\frac{M+1}{10}}\right).
\end{align*}
We are interested in the behavior of the polynomials $P^{[4]}_m$.

\begin{lemma} \label{PolyBound} For $j\in\{1,2,3\}$ we have for $1 \leq m \leq M$
	\begin{align*}
	\sup_{|x| \leq n^{\frac{1}{40}}} \left|\frac{P^{[j]}_m(x)}{n^{\frac{m}{10}}}\right| = O_M(1), \qquad n \rightarrow \infty.
	\end{align*}
	\end{lemma}
Setting $d_1:=\frac{1}{15}$, $d_2:=\frac{3}{80}$, and $d_3:=\frac{7}{80}$, the following lemma follows by a direct calculation.
\begin{lemma} \label{approxP4} We have for $|x| \leq n^{\frac{1}{40}}$
\begin{align*}
 \prod_{j=1}^{3} \lp 1+ \sum_{m=1}^{M} \frac{P^{[j]}_m(x)}{n^{\frac{m}{10}}}+O_{M}\lp n^{-d_j(M+1)} \rp \rp
&= \sum_{m=0}^M \frac{P^{[4]}_m(x)}{n^{\frac{m}{10}}} + \sum_{m=M+1}^{3M} \frac{P^{[4]}_{M,m}(x)}{n^{\frac{m}{10}}} + O_M\left( n^{-\frac{3(M+1)}{80}}\right) 
\\&= O_M(1).
\end{align*}
\end{lemma}

We finish this section with the following lemma which follows from Theorem \ref{Jsup}. 

\begin{lemma} \label{Jestimate} Let $\kappa > 0$ be sufficiently small. Then we have for all $z \in B_\kappa(1)$
	\begin{align*}
	\exp\left( \mathcal{E}_n(z)\right) = 1 + O\left( \left( \frac{B}{n^{\frac{1}{5}}}\right)^{n^{\frac{2}{55}}} \right),
	\end{align*}
	where $B > 0$ is some constant. 
\end{lemma}

\section{Wright's Circle Method and the proof of Theorem \ref{Main}}\label{sec:CircMeth}

\subsection{Overview of the strategy}

We write for $n \in \N_0$, using Cauchy's Theorem
\begin{equation*} \label{Rcircle}
r(n) = \frac{1}{2\pi i} \int_{C_n} \frac{G(q)}{q^{n+1}} dq,
\end{equation*}
where $C_n$ is a closed path inside the unit circle surrounding zero counterclockwise exactly once. To estimate $r(n)$ we use Wright's Circle Method \cite{Wright0} and split $C_n$ into two arcs: the major and the minor arc. The major arc $C^{\mathrm{maj}}_n$ is placed in a neighborhood of $q=1$, which is the dominant ``cusp" of $G(q)$. As we show in Propositions \ref{ExpandMain} and \ref{Minor} below the major arc produces the asymptotic main term and that the minor arc term is negligible.
In view of Proposition \ref{ErrorHolo} (1) and Theorem \ref{Asymptotic} it is natural to choose the major arcs in a way such that $z$ (given through $q = e^{-z}$) lies in a cone $\mathcal{C}_\delta$ for some $\delta$ independent from $n$. This gives us sufficiently good control over the error term $E(\eta; z)$. We use the line segment $z =  2X^2 n^{-\frac{3}{5}}(1 + iy)$, $-1 \leq y \leq 1$. Defining the minor arc as
\begin{align*}
C_n^{\min} := \left\{ q = e^{-2X^2 n^{-\frac{3}{5}} + iv} : 2X^2n^{-\frac35} < |v| \leq \pi\right\},
\end{align*}
we obtain
\begin{align} 
\label{rint} r(n) = \frac{1}{2\pi i} \int_{C^{\mathrm{maj}}_n} \frac{G(q)}{q^{n+1}}dq + \frac{1}{2\pi i} \int_{C^{\mathrm{min}}_n} \frac{G(q)}{q^{n+1}}dq.
\end{align}
The key idea is to approximate $\mathrm{Log}(G(e^{-z}))$ by Theorem \ref{Asymptotic}. However, since this series does not converge, we have to cut it of at a specific bound which we choose to be $k_n$. We abbreviate 
\begin{align*}
r_n(w) := H_n\left(iS_n - iw; n^{-\frac{1}{10}}\right) \exp(\mathcal{E}_n(w)), 
\end{align*}
where $S_n$ is the sequence of saddle points of Lemma \ref{SaddlePoint} and the functions $H_n$, and $\mathcal{E}_n$ are defined in \eqref{Hn} and \eqref{E}, respectively. We obtain, using Theorem \ref{Asymptotic} (choosing $k_n - \frac12 < \eta_n < k_n + \frac12$) and putting $a_n := 2X^2 n^{-\frac{3}{5}}$ 
\begin{align}
\nonumber \frac{1}{2\pi i} \int_{C^{\mathrm{maj}}_n} \frac{G(q)}{q^{n+1}}dq & = \frac{a_n}{2\pi i} \int_{1-i}^{1+i} G\left( e^{-a_nw}\right) e^{na_nw} dw\\
\label{MajorFinal} & = \frac{2 X^{\frac{4}{3}}}{in^{\frac{2}{5}} } \int_{1-i}^{1+i} \frac{\exp\left( n^{\frac{2}{5}}f_n(iS_n - iw) \right) r_n(w)}{w^{\frac13}} dw, 
\end{align}
where the function $f_n$, is defined in \eqref{Fn}. Note that $n^{\frac25}$, $na_n$, and $a_n^{-\frac23}$ have the same power in $n$. In Subsection \ref{expmint} we find an asymptotic expansion for the integral in \eqref{MajorFinal}. In Subsection \ref{marcest} we show that the minor arc integral in \eqref{rint} is negligible. Together, this proves Theorem \ref{Main}. 

\subsection{Modifying the path of integration and estimating the tails}

The goal of this subsection is to prove the following. 

\begin{proposition}\label{CurveMod}
	We have, as $n \rightarrow \infty$, 
	\begin{multline*}
	 \int_{1-i}^{1+i} \frac{\exp\left( n^{\frac{2}{5}}f_n(iS_n - iw) \right) r_n(w)}{w^{\frac13}} dw \\
	  = \int_{S_n-in^{-\frac{7}{40}}}^{S_n+in^{-\frac{7}{40}}} \frac{\exp\left( n^{\frac{2}{5}}f_n(iS_n - iw) \right) r_n(w)}{w^{\frac13}} dw + O\left( A(n) \exp\left(  -  \frac{5X^2}{3} n^{\frac{1}{20}}\right) \right).
	\end{multline*}
\end{proposition}

We need the following lemma.

\begin{lemma}\label{CompactLemma}
	The functions $w \mapsto w^{-\frac13}$ and $w \mapsto r_n(w)$  are uniformly bounded on $[(1-i)S_n,1-i]$ and $\{ w \in \IC : -S_n \leq \mathrm{Im}(w) \leq -n^{-\frac{7}{40}}, \mathrm{Re}(w) = S_n\}$ for $n$ sufficiently large. 
\end{lemma}

We are ready to prove Proposition \ref{CurveMod}. 

\begin{proof}[Proof of Proposition \ref{CurveMod}]
	Using the symmetry of the integrand under conjugation, we write
	\begin{multline}\label{newcurve1} 
	\int_{1-i}^{1+i} \frac{\exp\left( n^{\frac{2}{5}}f_n(iS_n - iw) \right) r_n(w)}{w^{\frac13}} dw =  \int_{S_n-in^{-\frac{7}{40}}}^{S_n+in^{-\frac{7}{40}}} \frac{\exp\left( n^{\frac{2}{5}}f_n(iS_n - iw) \right) r_n(w)}{w^{\frac13}} dw  \\
	+ O\left(\left(\int_{(1-i)S_n}^{S_n-in^{-\frac{7}{40}}}+\int_{1-i}^{(1-i)S_n}\right) \left| \frac{\exp\left( n^{\frac{2}{5}}f_n(iS_n - iw) \right) r_n(w)}{w^{\frac13}}  dw \right| \right),
	\end{multline}
	where all curves are chosen as straight lines.  We first observe with  \eqref{rEsti} for $n$ sufficiently large
	\begin{equation*}
		\sup_{w \in [(1-i)S_n,1-i]} \left|\exp\left(n^{\frac25} f_n(iS_n - iw)\right)\right| \le \exp\left(4.8X^2 n^{\frac25}\right).
	\end{equation*}
	We conclude with Lemma \ref{CompactLemma}
	\begin{equation}\label{FirstError}
	\int_{1-i}^{(1-i)S_n} \left| \frac{\exp\left( n^{\frac{2}{5}}f_n(iS_n - iw) \right) r_n(w)}{w^{\frac13}} dw\right|
	\ll \exp\left( 4.8 X^2 n^{\frac{2}{5}}\right) = \exp\left(\left(A_1 - 0.2X^2\right) n^{\frac{2}{5}}\right).
	\end{equation}
In the same way we obtain 
	\begin{equation*}
		\int_{S_n-iS_n}^{S_n-in^{-\frac{7}{40}}} \left| \frac{\exp\left( n^{\frac{2}{5}}f_n(iS_n - iw) \right) r_n(w)}{w^{\frac13}} dw \right| 
		\ll \sup_{\substack{-S_n \leq \mathrm{Im}(w) \leq -n^{-\frac{7}{40}} \\ \mathrm{Re}(w) = S_n}} \left|\exp\left(n^{\frac{2}{5}}f_n(iS_n - iw)\right)\right|.
	\end{equation*}
	Using Lemma \ref{Saddlesup}, \eqref{separatez4}, and \eqref{Sn}, we obtain
	\begin{align}\label{fn}
	\sup_{\substack{-S_n \leq \mathrm{Im}(w) \leq -n^{-\frac{7}{40}} \\ \mathrm{Re}(w) = S_n}}  \left|\exp\left(n^{\frac{2}{5}}f_n(iS_n - iw)\right)\right| \ll A(n)\exp\left( -  \frac{5X^2}{3} n^{\frac{1}{20}}\right).
	\end{align}
	With \eqref{newcurve1}, \eqref{FirstError}, and \eqref{fn},  this combines to 
	\begin{multline*}
	\int_{1-i}^{1+i} \frac{\exp\left( n^{\frac{2}{5}}f_n(iS_n - iw) \right) r_n(w)}{w^{\frac13}} dw = \int_{S_n-in^{-\frac{7}{40}}}^{S_n+in^{-\frac{7}{40}}} \frac{\exp\left( n^{\frac{2}{5}}f_n(iS_n - iw) \right) r_n(w)}{w^{\frac13}} dw \\
	+O\left( \exp\left(\left(A_1 - 0.2X^2\right) n^{\frac{2}{5}}\right) \right) + O\left( A(n)\exp\left(- \frac{5X^2}{3}n^{\frac{1}{20}}\right) \right).
	\end{multline*}
	The proposition now follows by comparing the error terms.
\end{proof}

\subsection{Expanding the main integral}\label{expmint}

In this section we expand the major arc integral \eqref{MajorFinal} to obtain a refined expression for $r(n)$. An important tool is the Saddle Point Method. We put 
\begin{align}
\label{Cm} C_m := 2X^{\frac{4}{3}} \exp(-A_5) \int_{-\infty}^{\infty} P^{[4]}_m(x) \exp\left( -\frac{5X^2x^2}{3}\right) dx,
\end{align}
where the polynomials $P^{[4]}_m$ are defined in Lemma \ref{approxP4}. We prove the following lemma.

\begin{proposition} \label{ExpandMain} Let $M \in \N$. We obtain, as $n \to \infty$,
	\begin{multline*}
	r(n)  = \frac{A(n)}{n^{\frac{3}{5}}} \left( \sum_{0 \leq m \leq \frac{3M}{8}} \frac{C_m}{n^{\frac{m}{10}}} + O_M\left( n^{-\frac{3(M+1)}{80}}\right)\right) + O\left(\int_{C^{\mathrm{min}}_n} \left| \frac{G(q)}{q^{n+1}}dq \right| \right) 
	\\
	+ O\left(A(n) \exp\left( - \frac{5X^2}{3} n^{\frac{1}{20}}\right)\right).
	\end{multline*}
	\end{proposition}

\begin{proof}
	 By \eqref{rint}, \eqref{MajorFinal}, and Lemma \ref{CurveMod}, we have, substituting $x = n^{\frac{1}{5}}(iS_n - iw)$  
\begin{align}\label{NewIntegral}
r(n) & = \frac{2 X^{\frac{4}{3}}}{n^{\frac{3}{5}} }  \int_{-n^{\frac{1}{40}}}^{n^{\frac{1}{40}}} \frac{\exp\left( n^{\frac{2}{5}}f_n\left(\frac{x}{n^{\frac{1}{5}}}\right) \right) H_n\left(\frac{x}{n^{\frac{1}{5}}}; n^{-\frac{1}{10}}\right)  \exp\left(\mathcal{E}_n\left(S_n + \frac{ix}{n^{\frac{1}{5}}}\right)\right)}{\left(S_n + \frac{ix}{n^{\frac{1}{5}}}\right)^{\frac13}} dx \\
\nonumber & \hspace{3cm}  + O\left(\int_{C^{\mathrm{min}}_n} \left| \frac{G(q)}{q^{n+1}}dq \right| \right)
+O\left( A(n)\exp\left( - \frac{5X^2}{3} n^{\frac{1}{20}}  \right)\right).
\end{align}
Together with Lemmas \ref{fexpand}, \ref{cuberootexpand}, \ref{Jestimate}, \ref{approxP4}, and Proposition \ref{HFinal} (note that $S_n + ixn^{-\frac{1}{5}} \in B_\kappa(1)$  for all $|x| \leq n^{\frac{1}{40}}$ for all $n$ sufficiently large, where $\kappa$ is chosen sufficiently small to satisfy Lemma \ref{Jestimate}) we find that the integral in \eqref{NewIntegral}, including the factor $2X^{\frac{4}{3}}n^{-\frac{3}{5}}$, equals 
\begin{multline}\label{OMForm}
\frac{2X^{\frac{4}{3}}\exp\left(-A_5\right) A(n)}{n^{\frac{3}{5}}} \int_{-n^{\frac{1}{40}}}^{n^{\frac{1}{40}}} \left(  \sum_{m=0}^M \frac{P^{[4]}(x)}{n^{\frac{m}{10}}} + \sum_{m=M+1}^{3M} \frac{P^{[4]}_{M,m}(x)}{n^{\frac{m}{10}}} + O_M\left( n^{-\frac{3(M+1)}{80}}\right) \right) 
\\ \times
 \left( 1 + O\left( \left( \frac{B}{n^{\frac{1}{5}}}\right)^{n^{\frac{2}{55}}} \right) \right) 
 \exp\left( -\frac{5X^2x^2}{3}\right) dx.
\end{multline}
By Lemma \ref{Gauss} we obtain
\begin{align}\nonumber
&\hspace{-.1cm}\frac{2X^{\frac{4}{3}}\exp\left(-A_5\right) A(n)}{n^{\frac{3}{5}}} \hspace{-.1cm} \int\limits_{-n^{\frac{1}{40}}}^{n^{\frac{1}{40}}} \hspace{-.2cm} \left(  \sum_{m=0}^M \frac{P^{[4]}_m(x)}{n^{\frac{m}{10}}} + \sum_{m=M+1}^{3M} \frac{P^{[4]}_{M,m}(x)}{n^{\frac{m}{10}}} + O_M\left( n^{-\frac{3(M+1)}{80}}\right) \right) \hspace{-.05cm} \exp\left( -\frac{5X^2x^2}{3}\right) dx \\
\label{ToConstants}
&\hspace{7.79cm}= \frac{A(n)}{n^{\frac{3}{5}}} \left( \sum_{0 \leq m \leq \frac{3M}{8}} \frac{C_m}{n^{\frac{m}{10}}} + O_M\left( n^{-\frac{3(M+1)}{80}}\right)\right),
\end{align}
where the $C_m$ are defined in \eqref{Cm}.
The claim can then be concluded using \eqref{NewIntegral}, \eqref{OMForm}, \eqref{ToConstants}, \eqref{Cm}, and Lemma \ref{approxP4}.
\end{proof}

\subsection{The minor arc estimate}\label{marcest}
Using Lemma \ref{Gestimate}, Proposition \ref{ErrorHolo} (1), and Theorem \ref{Asymptotic} yields the following estimate on the minor arc. 

\begin{proposition}[Minor arc estimate]\label{Minor}
	There exists $\delta > 0$ independent from $n$, such that 
	\begin{align*}
	\int_{C_n^{\mathrm{min}}} \left| \frac{G(q)}{q^{n+1}} dq \right| = O\left( n^{\frac{1}{5}}\exp\left( A_1n^{\frac{2}{5}} - \left(A_2 + \frac{\delta}{\sqrt{2}X}\right)n^{\frac{3}{10}} \right) \right), \qquad n \rightarrow \infty.
	\end{align*}
\end{proposition}

\subsection{Proof of Theorem \ref{Main}}

Theorem \ref{Main} now follows by a direct calculation using Propositions \ref{ExpandMain} and \ref{Minor}.

\subsection{The constants $C_j$}

We finally calculate the first constants  $C_0, C_1,$ and $C_2$ from Theorem \ref{Main}. Since we have $P^{[4]}_0 = 1$ by \eqref{approxP4}, by \eqref{Cm} 
\begin{align*}
C_0 = 2X^{\frac{4}{3}} \exp\left(-A_5\right) \int_{-\infty}^\infty \exp\left( - \frac{5X^2u^2}{3}\right) du = \frac{2\sqrt{3\pi} X^{\frac{1}{3}}}{\sqrt{5}} \exp\left(-A_5\right),
\end{align*}
where we use the well-known formula
\begin{align*}
\int_{-\infty}^\infty \exp\left( -bu^2\right) du = \sqrt{\frac{\pi}{b}}, \qquad b > 0.
\end{align*}
Note that this constant was already computed by Romik \cite{Ro}.

In the following, we provide the calculations to identify the next constants $C_1$ and $C_2$. The proof of Theorem \ref{Main} tells us that we require the first terms of the expansions in Lemmas \ref{fexpand}, \ref{cuberootexpand}, and Proposition \ref{HFinal} explicitly. We obtain with the help of a computer
\begin{align*}
P^{[4]}_1(x) & = -\frac{Y \left(35 x^2 X^2-6\right)}{120 X^3}-\frac{4959 Y^5}{102400000 X^{13}},\\
P^{[4]}_2(x) & = \frac{1}{27}  \left(40 x^2 X^2-9\right)i x +\frac{57 Y^6 \left(1015 x^2
	X^2-622\right)}{4096000000 X^{16}} +\frac{Y^2 \left(245 x^4 X^4-426 x^2
	X^2+36\right)}{5760 X^6}\\
&\hspace{10.25cm}+\frac{24591681 Y^{10}}{20971520000000000 X^{26}}.
\end{align*}
We use these polynomials, to calculate with \eqref{Cm}
\begin{align*}
\int\limits_{-\infty}^\infty P^{[4]}_1(x) \exp\left( - \frac{5X^2}{3} x^2 \right) dx &= - \sqrt{\frac{3 \pi }{5}}  \left(\frac{4959 Y^5}{102400000 X^{14}}-\frac{3 Y}{80 X^4}\right),\\
\int\limits_{-\infty}^\infty P^{[4]}_2(x) \exp\left( - \frac{5X^2}{3} x^2 \right) dx &= \sqrt{\frac{3 \pi }{5}} \left( \frac{24591681 Y^{10}}{20971520000000000 X^{27}}-\frac{7239 Y^6}{1638400000 X^{17}}-\frac{57 Y^2}{12800 X^7}\right).
\end{align*}
In particular we have 
\begin{align}
\label{C1} C_1 & = - 2X^{\frac43} \exp(-A_5)  \sqrt{\frac{3 \pi }{5}}  \left(\frac{4959 Y^5}{102400000 X^{14}}-\frac{3 Y}{80 X^4}\right), \\
\label{C2} C_2 & = 2X^{\frac43} \exp(-A_5) \sqrt{\frac{3 \pi }{5}} \left( \frac{24591681 Y^{10}}{20971520000000000 X^{27}}-\frac{7239 Y^6}{1638400000 X^{17}}-\frac{57 Y^2}{12800 X^7}\right).
\end{align}

\section{Open questions}
Although Theorem \ref{Main} provides an infinite number of terms in the asymptotic expression for $r(n)$, there are still interesting open questions. It is natural to ask about the nature of the constants $C_j$, for example, an explicit formula as rational functions in $X$, $Y$ and rational zeta and Gamma values (up to a constant $\exp(-A_5)$). 
Another question refers to the nature of the error term. We show that the error is smaller than any inverse polynomial times the main term $A(n)$, but it could be interesting to refine the suggested method.


\begin{thebibliography}{99}
	
	\bibitem{Andrews} G. Andrews, R. Askey, and R. Roy, {\it Special Functions}, Encyclopedia of Mathematics and its Applications {\bf 71}, Cambridge University Press, 2000.
	
	
	\bibitem{Apostol} T. Apostol, {\it Introduction to Analytic Number Theory}, Springer Science + Business Media, 1976.
	
	\bibitem{Brue} J. Brüdern, {\it Einführung in die analytische Zahlentheorie}, Springer, 1995. 
	
	
	%
	%
	
	
	
	
	\bibitem{Hall} B. Hall, {\it Lie Groups, Lie Algebras, and Representations: An Elementary Introduction.} Springer, 2004.
	
	\bibitem{HardyRama} G. Hardy and S. Ramanujan, {\it Asymptotic Formulæ in Combinatory Analysis}, Proceedings of the London Mathematical Society (2) {\bf 17} (1918), 75–115.
	
	
	
	\bibitem{Kaup} L. Kaup and B. Kaup, {\it Holomorphic Functions of Several Variables}, de Gruyter Studies in Mathematics {\bf 3}, 1983.
	
	\bibitem{Mat} K. Matsumoto, {\it On the analytic continuation of various multiple zeta- functions.} In: Number Theory for the Millenium II: Proceedings of the Millennial Conference on Number Theory (Urbana-Champaign, USA, 2000), 417–440. Eds. B. C. Berndt, N. Boston, H. Diamond, A. Hildebrand. AK Peters, 2002.
	
	
	
	\bibitem{Rademacher} H. Rademacher, {\it On the Partition Function $p(n)$}, Proceedings of the London Mathematical Society (2) {\bf 43} (1937), 241–254.
	
	\bibitem{Ro} D. Romik, {\it On the number of $n$-dimensional representations of $\SU(3)$, the Bernoulli numbers, and the Witten zeta function}, arXiv:1503.03776.
	
	
	
	
	\bibitem{Tenenbaum} G. Tenenbaum, {\it Introduction to Analytic and Probabilistic Number Theory}, Graduate Studies in Mathematics, American Mathematical Society \textbf{163}, Third Edition, 2008.  
	
	\bibitem{Wright0} E. Wright, {\it Asymptotic partition formulae II. Weighted partitions}, Proceedings of the London Mathematical Society {\bf 36} (1933), 117–141.
	



\end{thebibliography}
\end{document}